\documentclass[jmps]{apjrnl}
\usepackage{amssymb,amsfonts,amsmath}
{\theoremstyle{plain}%
  \newtheorem{theorem}{Theorem}
  \newtheorem{corollary}{Corollary}
  \newtheorem{proposition}{Proposition}
  \newtheorem{lemma}{Lemma}%
{\theoremstyle{remark}

}
{\theoremstyle{definition}

}
\setlength{\textwidth}{16.5cm}
\setlength{\textheight}{21.0cm}
\setlength{\topmargin}{-1cm}
\setlength{\leftmargin}{3cm}

\begin{document}

%% Author, fill in article title here
\title{On Successive Approximations \\ To The Choice Problem and Logic}

%% Fill in author list here
\author{Valeriy K. Bulitko
\\UTE Labs
\\4369 -- 31 Street
\\Edmonton, Alberta T6T 1C2
\\CANADA
\\E-mail: booly@shaw.ca
}

\date{January 20, 2002}
%please do not use \today, use actual date of version

\maketitle

\begin{abstract}
This paper studies the formation of logical operations from
pre-logical processes. We are concerned with the reasons for
certain mental processes taking form of logical reasoning and the
underlying drives for consolidation of logical operations in human
mind. Starting from Piaget's approach to Logic (Piaget, 1956) we
discuss whether the evolutionary adaptation can be such a driving
force and whether the limits of human mind can result in the
standard system of logical operations. The paper demonstrates that
the classical two-valued propositional logic can begin from a
method of successive approximations applied to a decision-making
problem within a framework of Subject-in-an-environment survival.
The presented results shed a new light on the known model of human
choice by Lefebvre (Lefebvre, 1991, 1995).
\end{abstract}


%%%% Authors begin text of article here %%%

\section{Introduction}

This paper studies the formation of logical operations from
pre-logical processes. We are concerned with the reasons for
certain mental processes taking form of logical reasoning and the
underlying drives for consolidation of logical operations in human
mind.

In order to enable a wider appreciation of the results, we opted
to present our model on the basis of Piaget's approach to Logic
(Piaget, 1956). Thus, even though this research was carried out
independently, the reader will find numerous references to
Piaget's concepts and definitions throughout the paper.

According to J.Piaget, logic is not an inherent form of thinking
and "logical operations result from coordinations of the actions
of combining, dissociating, ordering and setting up of
correspondences, which then acquire the form of reversible
systems" (Piaget, 1956, p. 13). In his works Piaget suggested a
representation of logical operations in terms of certain more
elementary ("pre-logical)") operations of the kinds listed above.

On the other hand, logical operations are the only means to carry
out logical inference. Therefore, it seems to be reasonable to
assume that the consolidation of logical operations in a
developing mind goes within a framework of some mental processes
which gradually take on the form of logical reasoning. What could
be the cause of such a transformation? How is it possible to
describe in strict mathematical terms such a "pre-logical" process
that comes to logical inference at the end?

In this paper, we attempt to demonstrate (i) that the
evolution-caused adaptation necessity could play the role of the
underlying driving force; (ii) that the processes of a successive
approach to the adaptation purpose can be behind logical
inference; and (iii) that then the limited nature of human mind
leads to the ordinary systems of logical operations.

The adaptive character of human intelligence is commonly accepted.
So it is natural to begin with the acceptance of an elementary
(i.e., as commonly applicable as possible) scheme of adaptation
problems and some fundamental solution processes for them.
%such that could model the logic (operations and processes) origin.
Such a scheme should be as elementary as to admit a formulation in
terms of the pre-logical operations and relations Piaget listed.

This research project studies one such scheme based on a model of
Subject's behavior choice (Bulitko, 2000). The model background
lies with two preference relations on the same set of Subject's
states. They are referred to as the {\em internal (subjective)
preference relation} and the {\em external (objective) preference
relation}.

The internal preference relation is defined by a partial order on
a set of inner states of the Subject. The external preference
relation on the state set is defined by a mapping from the state
set to another partially ordered set. The mapping will be called
{\em evaluation mapping (evaluation function)}.

We consider the process of getting to a better state with respect
to the external preference in a fashion that uses the external
preference relation as little as possible. On the contrary, the
Subject may use the internal preference without any restriction.

Several interpretations of this framework are possible. In our
primary interpretation the Subject must pay for its access to the
external preference whenever the access to the internal preference
is free. For example, we can think of the true external preference
as induced by a computationally expensive value function defined
over states while the subjective internal preference relation is
easy to compute. Thus, we strive to minimize the access to the
former while allowing unlimited access to the latter.

We use a successive approximation method as the core of our
solution approach. Humans are constantly involved in problem
solving that requires planning and execution of action {\em
sequences}. Therefore, it is reasonable to make the conjecture
that an apparatus enabling multi-step planning is among humans'
innate functions.

The successive approximation method proposed herein strongly
connects with Piaget's actions of "combining and dissociating". In
our model functions $\boxplus,\boxminus$ are the counterparts to
these actions. Since these pre-logical actions are ambiguously
described in (Piaget, 1956) we introduce our counterparts in an
axiomatic fashion.

The successive approximation method decomposes the evaluation
mapping into a superposition of functions $\boxplus,\boxminus$ as
well as "easier" functions (the details will be given below). Such
superpositions are called {\em approximating forms}.

We are now able to interpret Piaget's "actions of combining,
dissociating, ordering and setting up of correspondences" with
operations $\boxplus,\boxminus$, the two preferences, and the
evaluation mapping correspondingly.

In this paper, several approximating forms are developed for
different systems of axioms. Furthermore, two of them are relevant
to the hypothetical origin of logic we are discussing in this
research. These two are the cases in which the external and
internal preferences take certain {\em canonical forms}. Indeed,
every finite partially ordered set can be embedded into an
appropriate Boolean cube. Thus, it is always reasonable to replace
a given external preference with a two-element preference of the
"admissible - inadmissible" or "good-bad" kind.

Certain standard systems of logical operations obey these axioms.
Then the corresponding approximating forms translate into logical
formulae. Finally, the subject's action sequences map to logical
inferences. Thus, in our framework it is possible to explain the
transition from pre-logical to logical form of mental processes
via a canonical simplification of the internal and external
preferences.

Developing a comprehensive theory of algorithm design on the basis
of approximating forms is beyond the focus of this paper. Yet, we
pose an example of such an algorithm design in the final section
of the paper.

The approach proposed in this paper illuminates the well known
model of human choice by Lefebvre (Lefebvre 1991,1995) from a
different angle. In particular, the primary formula of that model
follows from the our approach under certain natural assumptions.

\section{Behavior Choice with Two Preferences}

In this section we consider the Subject faced with a choice of a
state among a set of states in the environment. Some of the states
may be better than others and some are incomparable. Subject's
objective is to reach a satisfactory terminal state.

A fundamental feature of many real-world behavior problems is the
difference between the evaluations of a state before and after the
state is arrived at. We will now attempt to formalize this
phenomenon by introducing two preference relations over the state
space.

One of the two relations will specify the "internal" system of
values based on the Subject's internal representation (or model)
of the world (including the Subject itself). We argue that the
internal preference relation is intrinsic to the Subject's
mentality since the Subject perceives the world in terms of this
relation. Thus, there are no restrictions on the usage of this
internal relation.

The other, "external", relation is based on effects of state
choices. Thus, the external preference reflects the actual nature
of the interaction between the environment and the Subject and,
generally speaking, only a part of the external preference is
available to the Subject. Naturally, the known portion of the
external preference relation includes the information the Subject
has uncovered so far in its exploration of the environment.

The external preference is objective and determines the Subject
rewards/penalties and ultimately its survival. No specific
limitations are imposed on the two preference relations making the
framework quite general. It is natural to pose an extremum (as
defined by the external preference) as a goal state for the
Subject. Indeed, maximum and the greatest elements can be
expressed in terms of preference relation by means of predicate
logic. It is clear reaching an extremum is a simpler problem than
reaching the optimal state.

Discrepancies between the internal and external preference
relations may cause problems for Subject. Note that there is a
cost associated with accessing the external preference relation.
It is not only the cost of accessing the information but also the
cost of changing the Subject's behavior patterns. (Here we
abstract from computing the actual cost values).

Thus informally, the problem studied herein is to find an extremum
of the external preference relation under certain given
restrictions on the information access to the external preference
and an unlimited access to the internal preference relation.

Given the restrictions and costs associated with accessing the
true (external) preference relation, the Subject strives to reach
its goal (i.e., to locate an extremum) using the internal
preference as much as possible. Naturally, in order for the
internal preference to be beneficial to the subject, it needs to
approximate the external preference. This interpretation of the
choice problem corresponds to a certain conservatism on the side
of the Subject when it is necessary to follow a certain external
pressure. Indeed, even if the Subject is aware of its incomplete
and/or incorrect representation it often might not be able to
correct it. Therefore, it will need to refine/reconstruct its
representation starting with whatever is available.

In order to address these problems, we consider a successive
approximation principle that will guide our further investigation.
Namely, in the following we will demonstrate that the problem is
decidable by some version of successive approximation method. The
underlying idea of the method is as follows. The Subject follows a
certain part of its internal preference as long as the preference
doesn't deviate significantly from the external one. Then on the
basis of accessible information on the external preference the
Subject reverses the corresponding part of the internal reference
and uses it to explore the environment further. The process then
repeats.

Thus, the Subject needs a {\em means} and a {\em scheme} to select
and manipulate corresponding parts of the preference relations.
The following section is devoted to a theory of such schema. We
believe the framework proposed below can be viewed as a possible
formalization of operations and relations listed by Piaget.

\section{An explication of the successive approximation method}

Let $S$ be a set of Subject states, $(M,\le_M),(L,\le_L)$ be
partially ordered sets of internal and external estimates
correspondingly. Let $\varphi:S\to M, \psi:S\to L$ be mappings
that link corresponding estimates to states. In this way it is
possible to set internal and external preferences on $S$.
Generally speaking, these preferences are pre-order relations on
the set (Birkhoff, 1967).

We simplify this description by introducing an order $\le_s$ on
$S$ through the mapping $\varphi$ and the poset $(M,\le_M)$ as
follows. Let us set $s\le_ss'\iff\varphi(s)<_M\varphi(s')\vee
s=s'$. This reformulation does not put any restrictions on
$\le_s$. Therefore, we can from now on consider the description
$\langle(M,\le_M),(L,\le_L), \psi:M\to L\rangle$. Furthermore,
$(M,\le_M)$ plays the role of $(S,\le_S)$ above and $\psi$ is
called the evaluation mapping (evaluation function).

It is worth noting that generally in each instance of the choice
problems the Subject gets the corresponding internal and external
preferences and the evaluation function. These three objects can
depend functionally on some parameters of the choice problem.

First, we consider the case of a single problem of choice. In a
section related to Lefebvre's model, we will consider a family of
choice problems.

If the evaluation function $\psi$ is a monotonic mapping (i.e.,
the condition $(\forall m_1,m_2\in
M)[m_1\le_Mm_2\Rightarrow\psi(m_1)\le_L\psi(m_2)$] is met) then
both preference relations $(M,\le_M)$ and $(L, \le_L)$ are
mutually compatible (concordant) and de facto the Subject may
follow its internal preference to reach the target state (that is,
a state with the maximum value).

Otherwise, it is natural to represent $\psi$ by a superposition of
monotonic evaluation mappings from $(M,\le_M)$ to $(L,\le_L)$ and
several connecting operations. We look for representations that
can be used for successive approximations. In finding a
representation of this kind that uses as few monotonic evaluation
mappings as possible, we attempt to use the external preference
relation as little as possible.

\subsection{Axiom system $\mathcal A$}

This section proposes collections of operations providing
representations of the aforementioned kind for any given
evaluation function $\psi$. These representations are called
"approximating forms".

Our first collection uses three operations: $\boxminus:L\times
L\to L,\boxplus:2^L\to L,\circledcirc:L\to L$. In our model, the
first two represent Piaget's operations of "dissociating" and
"combining" respectively. The third operation $\circledcirc$
represents the conception of the "null" element $o$ that we also
encounter in (Piaget, 1956). Following his theory, Piaget
developed a special algebra of numerous concrete operations.

We feel it is quite natural to define the sought model via an
appropriate axiomatic system. We start with system $\mathcal A$:
\begin{description}
\item[$\mathcal A_1$:]
       $(\forall S\subseteq M)(\exists\breve S\subseteq S)
       [(\forall s\in S)(\exists\breve s\in\breve S)
       [\breve s\le_Ms]\ \&\ (\forall\breve s,\breve s'\in\breve S)
       [\breve s\not<_M\breve s']]]$.
\item[$\mathcal A_2$:]
       $(\forall L',L''\subseteq L)(\forall x\in L')[(x\le_L\boxplus(L'))\&
       (L'\subseteq L''\Rightarrow\boxplus(L')\le_L\boxplus(L''))]$.
\item[$\mathcal A_3$:]
       $(\forall l,l'\in L)[\boxminus(l,\circledcirc(l))=l\ \&\
       (l\le_Ll'\Rightarrow\circledcirc(l)\le_L\circledcirc(l'))]$.
\item[$\mathcal A_4$:]
       $(\forall l,l'\in L)[l\le_Ll'\Rightarrow(\exists l''\in L)
       [\boxminus(l',l'')=l\ \&\ \circledcirc(l')\le_Ll'']]$.
\end{description}

Axiom $\mathcal A_1$ demands the internal preference to have no
infinite decreasing chains. The axiom is trivially true for finite
state set. It is clear that the restriction of finite state sets
is not overly constraining in practice. It is worth noting that
there is only one axiom relating to the internal preference.

Axiom $\mathcal A_2$ describes "combining" $\boxplus$ whereas
$\mathcal A_3,\mathcal A_4$ tie operation of "dissociating"
$\boxminus$ and operation $\circledcirc$ of coming to a "null".

Operation $\boxplus$ combines element set $L'\subseteq L$ into a
single element while respecting the monotonicity property. This
property is one of the main properties of set-theoretical
operation $\cup$. Thus, our definition preserves the primary
property of the concept of combination as used by Piaget.

Axiom $\mathcal A_4$ postulates the property of reversibility for
"dissociating". Element $l''$ such that $\boxminus(l,l'')=l'$
represents the "difference" between $l$ and $l'$ (again, this
preserves the flavor of Piaget's definition).

Axiom $\mathcal A_3$ fixes some sufficient properties of the
concept of "null". Note that many "null" elements may exist (but
not required to).

\subsection{Approximating forms}

For every function $\nu:M\to L$ we call set $\frak n(\nu)=
\{(m,m')|(m\le_Mm')\ \&\ (\nu(m)\not\le_L\nu(m'))\}$ {\sl
non-monotonicity domain of $\nu$}. If $\frak n(\nu)=\emptyset$
then $\nu$ is called {\sl monotonic function}. Also for every
poset $(R,\le_R)$ the standard mappings
$(\cdot)^{\vartriangle},(\cdot)^{\triangledown}:R\to2^R$ are
defined by $t^{\vartriangle}=\{t'\in R|t'\le_R t\},
t^{\triangledown}=\{t'\in R|t\le_R t'\}.$

\begin{theorem}\label{thm.thm1} Let for $(M,\le_M),(L,\le_L)$
all axioms of the system $\mathcal A$ be satisfied  and lengths of
all increasing chains in $(M,\le_M)$ do not exceed some integer
$D$. Then for every function $\psi:M\to L$ there exists a
representation $\psi=\boxminus(\varphi_1,
\boxminus(\varphi_2,\boxminus(\varphi_3,\dots)))$ such that all
$\varphi_i,i=1,2,3\dots,$ are monotonic functions from $(M,\le_M)$
to $(L,\le_L)$.

Furthermore, the number of occurrences of operation $\boxminus$ in
this representation does not exceed $D$.
\end{theorem}

\begin{proof} Let us reduce the problem for a given function
$\psi$ to the same problem for a simpler function $\psi_1$ such
that the following holds $\psi=\boxminus(\varphi_1,\psi_1)$ and
$\frak n(\psi_1) \subsetneqq\frak n(\psi)$.

First, we define $M_1=\{x\in M|\frak n(\psi)\cap(x^{\vartriangle}
\times x^{\vartriangle})\neq\emptyset\}$, $M^1=\overline{M_1}$ :
\begin{center}
$\varphi_1(x) =
\begin{cases}
\boxplus(\psi(x^{\vartriangle})),&  x\in M_1,\\
\psi(x),&   x\in M^1.
\end{cases}$
\end{center}
Then we set $\psi_1(x)$ to any such $z\in L$ that
$\boxminus(\varphi_1(x),z)=\psi(x)\ \&\
\circledcirc(\varphi_1(x))\le_Lz$ if $\varphi_1(x)\neq\psi(x)$.
Otherwise, we set $\psi_1(x)=\circledcirc(\psi(x))$.

Existence of element $z$ in the definition is guaranteed by axioms
$\mathcal A_3,\mathcal A_4$. Now, the equality
$\psi(x)=\boxminus(\varphi_1(x),\psi_1(x))$ holds due to the
definitions of $\varphi_1,\psi_1$.

Let us prove that function $\varphi_1:(M\le_M)\to(L,\le_L)$ is
monotonic.

First, $\varphi_1=\psi$ over $M^1$ and we may use the condition
$x,y\in M^1\& x\le_My\Rightarrow\psi(x)\le_L\psi(y)$. Indeed,
otherwise $\psi(x)\not\le_L\psi(y),x\le_My,\psi(x)\neq\psi(y)$
and, therefore, $y\in M_1\cap M^1$. However, $M^1\cap
M_1=\emptyset$ which leads to a contradiction.

Second, $\varphi_1$ maps $(M_1,\le_M)$ into $(L,\le_L)$
monotonically in accordance with $\mathcal A_2$.

Finally, let us consider the "mixed" case when $x\in M^1,y\in M_1$
and all elements of $M$ are comparable with respect to $\le_M$. It
is clear that $y\le_Mx$ is impossible since $z\in M_1\Rightarrow
z^{\triangledown}\subseteq M_1$ immediately follows from the
definition of $M_1$.

Thus, it remains to consider the possibility of $x\le_My$. In that
case $\varphi_1(y)=\boxplus(\psi(y^{\vartriangle}))\ge_L\psi(x)$
in accordance to $\mathcal A_2$. On the other hand,
$\psi(x)=\varphi_1(x)$ on $M^1$ follows from the definition of
$\varphi_1$. Hence, function $\varphi_1$ is monotonic.

We are now ready to prove the last assertion of the theorem. For
that it is sufficient to demonstrate the inclusion
$M^1\cup\breve{M_1}\subseteq M^2$. Here  $M^2,M_2$ are defined for
$\psi_1$ in the same way as $M^1,M_1$ were defined for $\psi$
above. $\breve{M_1}$ is the set of all minimal elements of set
$M_1$, see $\mathcal A_1$. Namely: $M^2=\overline{M_2}$ and
$M_2=\{x\in M|\frak n(\psi_1)\cap(x^{\vartriangle} \times
x^{\vartriangle})\neq\emptyset\}$.

From here we have $M_2\subseteq(M_1\setminus\breve{M_1})$ and
$\frak n(\psi_1)\subseteq\frak n(\psi)\setminus\breve M_1\times
M_1$. Thus, sequence $M_1\supsetneqq M_2\supsetneqq
M_3\supsetneqq\dots$ ends on a step with the number can not be
greater than the highest of the lengths of the increasing chains
in poset $(M,\le_M)$. Indeed, since $\breve M_2\subseteq
M_1\setminus\breve M_1$ then in accordance with $\mathcal A_1$ for
every element $y\in\breve M_2$ there exists some $x\in\breve M_1$
such that $x<_Ly$. Therefore, one can choose an increasing chain
of representatives of sets $\breve M_1,\breve M_2,\breve
M_3,\dots$ which are mutually disjoint sets.

We will now prove that $M^1\cup\breve{M_1}\subseteq M^2$. First,
$\varphi_1(x)=\psi(x)$ holds for every $x\in M^1$. From here
$\psi_1(x)=\circledcirc(\psi(x))$. However, mapping $\psi_1$ is
monotonic on $M^1$ due to $\mathcal A_3$ and since $\psi$ is
monotonic on $M^1$. So $(M^1\times M^1)\cap\frak
n(\psi_1)=\emptyset$ and therefore $M^1\subseteq M^2$.

Further, let $x,y\in M^1\cup\breve{M_1}$ and $x\le_M y$. Then we
can show that $\psi_1(x)\le_L\psi_1(y)$. Indeed, the case $x,y\in
M^1$ was considered above. The case $x,y\in\breve{M_1}$ is
impossible since all elements of $\breve{M_1}$ are incomparable by
the definition. Above, we saw that $x\in M_1\ \&\
x\le_My\Rightarrow y\in M_1$. Besides $M^1\cap M_1=\emptyset$.
Therefore, $x\in M^1,y\in\breve M_1$ is the only case remaining to
consider. By definition $\psi_1(x)=\circledcirc(\psi(x))$ and
relation $\boxminus(\varphi_1(y),\psi_1(y))=\psi(y)$ holds.
Moreover, $\psi(y)<_L\varphi(y)$. In accordance with $\mathcal
A_4$ we have $\circledcirc(\varphi_1(y))\le_y\psi_1(y)$. Hence,
$\psi_1(x)\le_L\psi_1(y)$ takes place since $\circledcirc$ is a
monotonic operation due to $\mathcal A_3$ and
$\psi(z)\le_L\varphi_1(z),z\in M$ in accordance to $\mathcal A_2$
and the construction. \end{proof}

Let us denote by $\mathcal M$ the class of all monotonic mappings
from $(M,\linebreak\le_M)$ to $(L,\le_L)$. Also let $S(\varphi)=
\{x|\varphi(x)>_L\circledcirc(x)\},\varphi\in\mathcal M$.

\begin{corollary}\label{thm.cor1} Under the conditions of
theorem \ref{thm.thm1} for every $\psi:M\to L$ there exists a
substitution $p:\{z_1,\dots,z_{D+1}\}\to\mathcal M$ such that
$S(p(z_{n+1}))\subseteq S(p(z_{n})), n=\overline{1,D}$, and
\begin{center}
$\psi=\bold{Sb}^{z_1\ \ \dots\ \ z_{D+1}}_{p(z_1)\dots p(z_{D+1})}
\boxminus(z_1,\boxminus(z_2,\boxminus(\dots\boxminus(z_D,z_{D+1})\dots)))$.
\end{center}
\end{corollary}

\begin{proof} Let us fix a formula $\Phi(z_1,\dots,z_{D+1})=
\boxminus(z_1,\boxminus(z_2,(\dots
\boxminus(z_D,z_{D+1})\dots)))$ and consider substitutions of
monotonic functions instead of variables $z_1,\dots,z_{D+1}$  when
their results are determined.

According to theorem \ref{thm.thm1}, for every $\psi:M\to L$ there
exists representation
$$
\psi=\boxminus(\varphi_1,\boxminus(\varphi_2,
\boxminus(\varphi_3,\dots)))
$$
where all $\varphi_i,i=1,2,3\dots$
are monotonic mappings from $(M,\le_M)$ into $(L,\le_L)$.
Condition $S(\varphi_{n+1})\subseteq S(\varphi_n)$ follows the
construction of functions $\varphi_i$ made in the proof of theorem
1. The number of occurrences of operation $\boxminus$ in this
representation does not exceed $D$.

Once representation:
$$
\psi=\boxminus(\varphi_1,\boxminus(\varphi_2,\boxminus(\varphi_3,
\dots\boxminus(\varphi_{k},\varphi_{k+1})\dots)))
$$
with $k<D$ is obtained for mapping $\psi$, one can always continue
the expression on the right side of the representation until
$k=D$. For that it is sufficient to set
$$
\varphi_{i}(x)=\circledcirc(\varphi_{i-1}(x)),i=\overline{k+2,D+1}.
$$
In accordance with axiom $\mathcal A_3$ the obtained functions are
monotonic and
$$
\psi=\boxminus(\varphi_1,\boxminus(\varphi_2,\boxminus(\varphi_3,
\dots\boxminus(\varphi_{D},\varphi_{D+1})\dots))).
$$
\end{proof}

For any given $D$ the corollary states existence of the universal
formula
$$
\boxminus(z_1,\boxminus(z_2,
\boxminus(\dots\boxminus(z_D,z_{D+1})\dots)))
$$
which describes a structure of the representations. However, the
cost of the universality lies with the fact that the length of the
representation in theorem \ref{thm.thm1} can be essentially lower
than the lengths of the representations suggested by the
corollary.

In order to apply the theory developed in the last section, we
need to specialize monotonic functions in approximating forms.
Thus, we define the functions via the following auxiliary
construction.

Let us denote by $R^{\bot}$ the set of all minimal elements of
$(R,\le_R)$.
On the basis of axiom $\mathcal A_1$ let us split set $M$:\\
$M_1=M^{\bot}$;\\
$M_{n+1}=(M\setminus\underset{j\le n}{\cup}M_j)^{\bot}$.\\
Any two elements of $M_j$ are incomparable in $(M,\le_M)$ for any
$j$.

We denote by $\theta$-function of rank $i$ any monotonic mapping
$\theta:M\to L$ such that $\theta(x)=\circledcirc(x)$ for all
$x\in\underset{j<i}{\cup}M_j$ as well as
$\theta(x)\in\max(L,\le_L)$ for all $x\in\underset{j>i}{\cup}M_j$.
The rank of a given function $\theta$ is denoted as
$\rho(\theta)$. Let $\Theta$ be the class of all
$\theta$-functions.

\begin{theorem}\label{thm.thm2}
Let conditions of theorem 1 be fulfilled, $D$ be the exact upper
bound of the lengths of the increasing chains in $(M,\le_M)$ and
$(L,\le_L)$ contain its greatest element. Then for any function
$\psi:M\to L$ there exists a substitution
$p:\{z_1,\dots,z_{D+1}\}\to\Theta$, such that
$\rho(p(z_i))=i,i=\overline{1,D+1}$ and $\psi=\bold{Sb}^{z_1\ \
\dots\ \ z_{D+1}}_{p(z_1)\dots p(z_{D+1})}
\boxminus(z_1,\boxminus(z_2,\boxminus(\dots\boxminus(z_D,z_{D+1})\dots)))$.
\end{theorem}

\begin{proof} We denote by $\gamma$ the greatest element of $(L,\le_L)$
and use induction on $D$. In the case of $D=0$ the statement is
obvious since there are no restrictions on $\theta$-functions.
Therefore, $\psi\in\Theta$.

{\sl Induction step:} Let us define $\theta_1$ of rank 1 in the
following manner:
\begin{center}
$\theta_1(x) =
\begin{cases}
\psi(x),&  \text{ if }x\in M_1,\\
\gamma,&   \text{ else}.
\end{cases}$
\end{center}
Then we may state $\psi=\boxminus(\theta_1,\psi_1)$ where for
$\psi_1$ we have $\psi_1(x)=\circledcirc(x)$ if $x\in M_1$ else
$\psi_1(x)$ satisfies $\boxminus(\gamma,\psi_1(x))=\psi(x)$.

It remains to obtain the desirable representation of $\psi_1$ on
set $M\setminus M_1$ with the partial order $\le_{m'}$ induced by
$\le_M$. Since the length of the longest increasing chain in
$(M\setminus M_1,\le_{m'})$ is $D-1$ we may use the induction
supposition. \end{proof}

\subsection{Axiom system $\mathcal B$.}

Let us define binary operations $\boxminus,\uplus:L\times L\to L$
and a unary operation $\circledcirc:L\to L$ in such a way that the
system $\mathcal B=\{\mathcal B_1,\dots,\mathcal B_4\}$ of axioms
takes place. Here $\mathcal B_i$ coincides with $\mathcal A_i$ for
$i=1,3,4$. Also:
\begin{description}
\item[$\mathcal B_2$:]
$(\forall x,y\in L)[x\le_L\uplus(x,y)\ \&\ y\le_L\uplus(x,y)]$.
\end{description}

\begin{theorem}\label{thm.thm3} Let for $(M,\le_M),(L,\le_L)$
all axioms of the system $\mathcal B$ be satisfied, lengths of all
increasing chains in $(M,\le_M)$ do not exceed some integer $D$,
and every increasing chain in $(L,\le_L)$ be a finite one. Then
for every $\psi:M\to L$ there exists a representation
$\psi=\boxminus(\varphi_1,
\boxminus(\varphi_2,\boxminus(\varphi_3,\dots)))$ where all
$\varphi_i,i=1,2,3\dots$ are monotonic  functions from $(M,\le_M)$
to $(L,\le_L)$.

The number of occurrences of the operation $\boxminus$ in this
representation does not exceed $D$. \end{theorem}

\begin{proof} First, in the case when $(\forall x\in
M)[|x^{\vartriangle}|<\infty]$ is true we can prove this theorem
using theorem 1. For that we will only need to note that in this
case it is possible to replace $\boxplus(\psi(x^{\vartriangle}))$
with any expression of the kind $\uplus(\psi(z_1),\uplus(\dots
\uplus(\psi(z_{n-1}),\psi(z_{n}))\dots))$. Here $z_1,\dots,z_n$ is
an enumeration of the finite set $x^{\vartriangle}$. Indeed, in
the proof of theorem \ref{thm.thm1} we used axiom $\mathcal A_2$
only for subsets of $L$ of the form $\psi(x^{\vartriangle})$.
Thus, it is sufficient to check that axiom $\mathcal A_2$ is
respected for sets of the kind $\psi(x^{\vartriangle})$. This
check is trivial on the basis of axiom $\mathcal B_2$ for
operation $\uplus$.

Otherwise, when there are infinite sets $x^{\vartriangle}$ we can
make use of the same scheme for the operation $\boxplus$ basing on
the condition of finiteness of increasing chains in $(L,\le_L)$.
For that let us enumerate elements $z_1,z_2,\dots,z_n,\dots$ of
set $x^{\vartriangle}$ for a given
$x\in M$. Simultaneously we compute a series of expressions:\\
$\uplus(\psi(z_1),\psi(z_2))$,\\
$ \uplus(\uplus(\psi(z_1),\psi(z_2)),\psi(z_3))$,\\
$\dots$\\
$\uplus(\uplus(\dots\uplus(\uplus(\psi(z_1),\psi(z_2)),
\psi(z_2))\dots),\psi(z_n))$,\\
$\dots$\ .

By axiom $\mathcal B_2$ the values of these expressions are
comparable and do not decrease in $(L,\le_L)$. In view of the
finiteness supposition for increasing chains in $(L,\le_L)$ the
sequence of computed values becomes stable from a certain element.
We set $\varphi_1(x)$ to this final value.

Thus, $\varphi_1$ is a monotonic mapping and
$\psi(x)\le_L\varphi_1(x),x\in M$. The last part of the proof is
analogous to the corresponding part of theorem \ref{thm.thm1}.
\end{proof}

\subsection{Dual problem}

Following the case of maximization investigated above, we now
consider the minimization problem whose formulation can be
obtained from the previous case simply via replacing the word
"maximization" with the word "minimization".

There is an easy reduction of the minimization problem to the
maximization. To do the reduction we first replace $\le$ with
$\ge$ in both preferences. Secondly, we replace functions
$\varphi_i$ with dual ones that remain to be non-decreasing
monotonic and then increasing chains with decreasing ones, maximal
elements with minimal ones, etc. This way we arrive at a new set
of axioms for operations denoted by an asterisk. System $\mathcal
A$ is replaced with $\mathcal A^{\star}$:
\begin{description}
\item[$\mathcal A^{\star}_1$:]
       $(\forall S\subseteq M)(\exists\breve S\subseteq S)
       [(\forall s\in S)(\exists\breve s\in\breve S)
       [\breve s\ge_Ms]\ \&\ (\forall\breve s,\breve s'\in\breve S)
       [\breve s\not<_M\breve s']]]$.
\item[$\mathcal A^{\star}_2$:]
       $(\forall L',L''\subseteq L)(\forall x\in L')[(x\ge_L
       \boxplus^{\star}(L'))\&
       (L'\subseteq L''\Rightarrow\boxplus^{\star}(L')\ge_L
       \boxplus^{\star}(L''))]$.
\item[$\mathcal A^{\star}_3$:]
       $(\forall l,l'\in L)[\boxminus^{\star}(\circledcirc^{\star}(l),l)=l
       \ \&\
       (l\ge_Ll'\Rightarrow\circledcirc^{\star}(l)\ge_L
       \circledcirc^{\star}(l'))]$.
\item[$\mathcal A^{\star}_4$:]
       $(\forall l,l'\in L)[l\ge_Ll'\Rightarrow(\exists l''\in L)
       [\boxminus^{\star}(l'',l')=l\ \&\
       \circledcirc^{\star}(l')\ge_Ll'']]$.
\end{description}
Then the following theorem, that is a dual for theorem
\ref{thm.thm1}, can be proven:

\begin{theorem}\label{thm.thm1*}
Let for posets $(M,\le_M),(L,\le_L)$ all axioms of the system
$\mathcal A^{\star}$ be respected and the lengths of all
decreasing chains in $(M,\le_M)$ do not exceed a certain integer
$D$. Then for every function $\psi:M\to L$ there exists
representation $\psi=\boxminus^{\star}(\varphi_1,
\boxminus^{\star}(\varphi_2,\boxminus^{\star}(\varphi_3,\dots)))$
where all $\varphi_i,i=1,2,\dots,$ are monotonic functions from
$(M,\le_M)$ to $(L,\le_L)$.

The number of occurrences of the operation $\boxminus^{\star}$ in
this representation does not exceed $D$. \end{theorem}

In the same manner we can formulate system $\mathcal B^{\star}$ of
axioms:
\begin{description}
\item[$\mathcal B^{\star}_1$:]
       $(\forall S\subseteq M)(\exists\breve S\subseteq S)
       [(\forall s\in S)(\exists\breve s\in\breve S)
       [\breve s\ge_Ms]\ \&\ (\forall\breve s,\breve s'\in\breve S)
       [\breve s\not<_M\breve s']]]$.
\item[$\mathcal B^{\star}_2$:]
       $(\forall x,y\in L)[x,y\ge_L\uplus^{\star}(x,y)]$.
\item[$\mathcal B^{\star}_3$:]
       $(\forall l,l'\in L)[\boxminus^{\star}(\circledcirc^{\star}(l),l)=l
       \ \&\
       (l\ge_Ll'\Rightarrow\circledcirc^{\star}(l)\ge_L
       \circledcirc^{\star}(l'))]$.
\item[$\mathcal B^{\star}_4$:]
       $(\forall l,l'\in L)[l\ge_Ll'\Rightarrow(\exists l''\in L)
       [\boxminus^{\star}(l'',l')=l\ \&\
       \circledcirc^{\star}(l')\ge_Ll'']]$.
\end{description}
Then a dual to theorem 3 holds:

\begin{theorem}\label{thm.thm3*} Let for posets $(M,\le_M),(L,\le_L)$
all axioms of the system $\mathcal B^{\star}$ be respected, the
lengths of all decreasing chains in $(M,\le_M)$ do not exceed some
integer $D$, and every decreasing chain in $(L,\le_L)$ be a finite
one. Then for every function $\psi:M\to L$ there exists
representation $\psi=\boxminus^{\star}(\varphi_1,
\boxminus^{\star}(\varphi_2,\boxminus^{\star}(\varphi_3,\dots)))$
where all $\varphi_i,i=1,2,\dots,$ are monotonic functions from
$(M,\le_M)$ to $(L,\le_L)$.

The number of occurrences of the operation $\boxminus^{\star}$ in
this representation does not exceed $D$.
\end{theorem}

It is said that $(L,\le_L)$ admits a dual isomorphism $\eta$ if
$\eta$ is an one-to-one mapping of $L$ onto $L$ such that $\forall
l,l'[l\le_Ll'\iff \eta(l')\le_L\eta(l)]$ holds.

If the external preference $(L,\le_L)$ admits a dual isomorphism
$\eta$ then the following identities hold:
$$
\chi^{\star}=\eta^{-1}\circ\chi\circ\eta,\chi\in\{\boxplus,\boxminus,
\circledcirc\}
$$
where $\circ$ denotes the composition of functions.

Thus generally speaking, we obtain new operation systems and new
representations that we refer to as approximating forms.

\section{A Possible Origin of Logic}

In this section we include a complexity notion into our
considerations. First of all, the Subject might reduce the
external preferences to the simplest kind such as
"acceptable-unacceptable" or "good-bad", etc. So in this case we
can set $L=\{0,1\},\le_L=\{(0,0),(0,1),(1,1)\}$.

Now it is natural to use the simplest collection of operations. As
well known, poset $(\mathcal B^n,\preccurlyeq)$ is a self-dual
poset for any  $n$. In particular, given aforementioned
$(L,\le_L),(n=1),$ we have $\eta(0)=1,\eta(1)=0$ with identity
$\eta=\eta^{-1}$. Theorems 3 and 3$^{\star}$ offer two-argument
operations $\uplus$ and $\uplus^{\star}$ correspondingly (unlike
many-place operations $\boxplus,\boxplus^{\star}$ from theorems
1,1$^{\star}$) for this case. Both representations introduced in
theorems 3, 3$^{\star}$ holds and $\eta(l)=\neg l$ is true.

\begin{lemma}\label{thm.lem1}  Let
$L=\{0,1\}$ and $\le_L=\{(0,0),(0,1),(1,1)\}$. Then:
\begin{enumerate}

\item Operation $\lambda x,y[\neg x\ \&\ y]$ as $\boxminus$,
operation $\Bbb{\bold0}:\{0,1\}\to\{0\}$ as
$\circledcirc^{\star}$, and operation $\vee$ as $\uplus$ obey the
axiom set $\mathcal B$.

\item Operation $\lambda x,y[y\ \to\ x]$ as $\boxminus^{\star}$,
operation $\Bbb{\bold1}:\mathcal B^n\to\{1\}$ as
$\circledcirc^{\star}$, and operation $\&$ as $\uplus^{\star}$
obey the axiom set $\mathcal B^{\star}$.

\item There exists only one boolean interpretation
of the operations $\boxminus,\boxminus^{\star}$.

\end{enumerate}
\end{lemma}

The lemma can be proved via a routine check of the axiom systems.

Let us recall that $x\to^{\star}y=\neg(y\to x)=\neg x\ \&\ y$ (see
for example (Kleene, 1967)). Henceforth, we refer to the
approximating forms constructed with operations
$\to^{\star},\vee,0$ or with $\to,\&,1$ as {\em boolean
approximating forms}.

The important question here is why natural human languages do not
contain any connective that represents operation $\to^{\star}$ (in
the way like the connective "and" represents $\&$, for example). A
possible answer is offered below.

Dual approximating forms of theorem 3$^{\star}$ begin with a given
function $\psi$ and approximate it by means of successive
simplifications: $\psi_i=
\boxminus^{\star}(\varphi_{i+1},\psi_{i+1})$, where $\psi_0=\psi$
and $i$ runs integers $0,1,2,\dots,$ while $\psi_i$ is not a
monotonic function (i.e., not an "easy" one). Taking in account
the meaning of $\boxminus^{\star}(x,y)$ is $y\to x$ we get
$\psi=\psi_1\to\varphi_1=(\psi_2\to\varphi_2)\to\varphi_1,\dots$.

Since $\varphi_1\to\psi$ follows from $\psi=\psi_1\to\varphi_1$ we
can think that the transition from $\psi$ to $\varphi_1$ means the
transition from the general notion $\psi$ to the specific notion
$\varphi_1$. On the contrary, in the dual case we have only
$\overline{\psi^{\star}\to\varphi_1^{\star}}$. Taking into account
the fact that a developing mind forms classes from specific
examples we see a support to the claim that the first transition
(from {\sl general} to {\sl specific}) is easier to implement.

Until now we have not assumed any properties about preference
$(M,\le_M)$. The second step of the simplification process is an
isotonic embedding a given finite internal preference $(M,\le_M)$
into an appropriate Boolean cube $\mathcal B^n$ where $n$ is the
dimension of the cube. This step is always possible for finite
preferences (Birkhoff, 1967). Therefore, let the given inner
preference be $(\mathcal B^n,\preccurlyeq)$. Then direct
corollaries of the lemma \ref{thm.lem1} and theorems above are as
follows:

\begin{corollary}\label{them.cor2}  Every classical logic function can be
represented by a boolean approximating form. \end{corollary}

\begin{corollary}\label{thm.cor3} Every $n$-argument logical (boolean)
function $f$ can be represented by an implicative normal form of
the kind  $f=P_k\to P_{k-1}\to\dots\to P_0$, where $k\le n$, and
$P_i,i=\overline{0,k},$ are monotonic boolean function.
\end{corollary}

Therefore, one can consider the classical two-valued propositional
logic merely as an application of the above-mentioned principle of
successive approximations to the problem of decision-making within
the Subject-environment survival framework. Thus, this viewpoint
suggests a way for the classical propositional logic to develop
from the survival problem. It is also important that this
hypothetical origin of logic appears quite natural.

\section{Application to one model by Lefebvre}

Lefebvre proposed (Lefebvre, 1991,1995) a model of Subject facing
a choice among a set of alternatives. In the model the Subject is
represented by function $X_1=f(x_1,x_2,x_3)$ where
$X_1,x_1,x_2,x_3$ run over $[0,1]$. The value of $X_1$ is
interpreted as "the readiness to choose a positive pole"
(Lefebvre, 1991) with probability $X_1$, and the value of $x_3$ as
the Subject's plan or intention to choose a positive pole with
probability $x_3$. Variables $x_1$ and $x_2$ represent the world
influence on the subject.

Furthermore, function $f$  is required to obey the following
axioms:
\begin{description}
\item[$\mathcal L_1$:]
$(\forall x_3\in[0,1])[f(0,0,x_3)=x_3]$\ ("the axiom of free
choice");
\item[$\mathcal L_2$:]
$(\forall x_3\in[0,1])[f(0,1,x_3)=0]$\ ("the axiom of credulity");
\item[$\mathcal L_3$:]
$(\forall x_2,x_3\in[0,1])[f(1,x_2,x_3)=1]$\ ("the axiom of
non-evil-inclinations");
\item[$\mathcal L_4$:]
$(\forall i,j,k)[\{i,j,k\}=\{1,2,3\})\Rightarrow (\forall
x_j,x_k\in[0,1])(\exists c,c'\in\Bbb R) (\forall
x_i\in[0,1])[f(x_1,$ $x_2,x_3)=cx_i+c']]$ ("the postulate of
simplicity").
\end{description}

\noindent Through this model Lefebvre gave explanations of several
psychological experiments putting it in the spotlight (e.g., see
bibliography in (Lefebvre, 1995)). It is, however, worth-while to
ponder if the model is mainly a compact representation (i.e., a
"roll-up") of certain empirical data or wether it describes a
fundamental structure governing human behavior.

In order to substantiate his model Lefebvre used various arguments
including the well known "antrophic principle" (Lefebvre, 1995).
In addition to our previous comments (Bulitko, 1997), in the
following we present an alternative justification to Lefebvre's
model rooted in the theory of the approximating forms presented in
the prior sections.

First, we show a reduction of the general case to the boolean
case. Second, we demonstrate that the system of the first three
axioms by Lefebvre can be replaced with a postulate of special
poset $(M,\le_M)$ and a special algorithm computing a decision
(choice). Namely, the poset can be chosen in the form of a linear
ordered three-element set. We furthermore suggest a natural
interpretation of such poset $(M,\le_M)$ and the algorithm.

\subsection{Lefebvre's ensembles}

It is easy to check that in the boolean case
($X_1,x_1,x_2,x_3\in\{0,1\}$) the axioms $\mathcal L_1-\mathcal
L_3$ completely define $f$. Namely, $f(x_1,x_2,x_3)=(x_3\to
x_2)\to x_1$. The "postulate of simplicity" $\mathcal L_4$ sets
$f$ on the interior of the three-dimensional cube $[0,1]^3$ in the
real-valued case.

Let us consider a set $Q$ of Subjects $s_i$  each being described
by a probabilistic collection $\tilde{\alpha}_i$ of values of the
boolean variables $(n_1,n_2,n_3)$. Let us assume that the
probability of encountering a Subject with a collection
$\tilde{\alpha}$ of the variable values in $Q$ is equal to
$p_{\tilde{\alpha}}$.

If behavior $z_i$ of each $s_i\in Q$ is described by the function
$(n_3\to n_2)\to n_1$ then we refer to $Q$ as {\em Lefebvre's
ensemble ($L$-ensemble or {\rm simply} ensemble) $\langle
Q,P\rangle$ with characteristic $P=(p_0,\dots,p_7)$}. We call
elements of the $L$-ensemble $L$-Subjects. Here $p_k$ denotes
$p_{\tilde{\alpha}}$ and  $k$ is the decimal representation of the
binary sequence $\tilde{\alpha}$.

Ensemble $\langle Q,P\rangle$ averaging boolean variables
$n_1,n_2,n_3,z_i$ yields real numbers $x_1,x_2,x_3,z\in[0,1]$.
Given the truth table of the boolean function $n_3\to n_2\to n_1$
elementary probabilistic considerations lead to the following
equalities:
\begin{eqnarray}
1=\overset{7}{\underset{k=0}{\Sigma}}p_k,\\
x_1=p_4+p_5+p_6+p_7,\\
x_2=p_2+p_3+p_6+p_7,\\
x_3=p_1+p_3+p_5+p_7,\\
z=p_1+p_4+p_5+p_6+p_7.
\end{eqnarray}
It is therefore reasonable to inquire which $L$-ensembles $\langle
Q,P\rangle$ values of $x_1,x_2,x_3,z$ satisfy Lefebvre's equation
$z=x_1+(1-x_1-x_2+x_2x_3)x_3$.

The following examples show that, generally speaking, $z\neq
f(x_1,x_2,x_3)$. Indeed, let us set
$p_1=p_2=p_3=p_4=p_5=p_6=p_7=0,1$. Then $x_1=x_2=x_3=0.4$ and
$f(x_1,x_2,x_3)=0.544$. However, the ensemble average $z$ equals
$0.5$. Interestingly enough, the difference can be quite
substantial as the following example demonstrates. Namely,
$p_0=p_1=p_2=p_4=p_6=p_7=0,p_3=p_5=0.5$ correspond to
$x_1=x_2=0.5,x_3=1$. Then $z=0.5$ but $f(0.5,0.5,1)=0.75$. Thus,
the error can reach at least 30\%.

On the other hand, equality $z=f(x_1,x_2,x_3)$ is met for all
possible (i.e., obeying equations (1)-(4)) characteristics $P$
when
$(x_1,x_2,x_3)\in\{(x_1,x_2,x_3)|x_1=1\}\cup\{(x_1,x_2,%\linebreak
x_3)|x_2=0\}\cup\{(x_1,x_2,x_3)|x_2=1\}\cup\{(x_1,x_2,x_3)|x_3=0\}$.

\begin{proposition}\label{thm.pro1}
For every collection $x_1,x_2,x_3\in[0,1]$ there exists an
$L$-ensemble $\langle Q,P(x_1,x_2,x_3)\rangle$ with characteristic
$P(x_1,x_2,x_3)$ such that $z=f(x_1,x_2,x_3)$.
\end{proposition}

\begin{proof} Let us consider three independent boolean
random variables $\zeta,\eta,\theta:\Bbb N\to\{0,1\}$ with the
mean values $x_1,x_2,x_3$ correspondingly. Then random variable
$(\zeta,\eta,\theta):\Bbb N\to\{0,1\}^3$ runs over the desired
ensemble $\langle Q,P(x_1,x_2,x_3)\rangle$. For the $i$-th
component of the characteristic $p_i(x_1,x_2,x_3)=
\underset{j=1,2,3}{\Pi}(1-\sigma_j+(-1)^{1-\sigma_j}x_j)$ holds
where $\sigma_j\in\{0,1\},j=\overline{1,3},$ and
$i=\underset{j=1,2,3}{\Sigma}\sigma_j2^{3-j}$. A simple
verification shows that the relations (1)-(4) are fulfilled and if
$z$ satisfies (5), then $z=f(x_1,x_2,x_3)$. \end{proof}

We call the ensembles described in this proposition {\it pure
Lefebvre's ensembles ($PL$-ensembles)}. Thus, a $PL$-ensemble is a
collection of $L$-Subjects with random parameters $(n_1,n_2,n_3)$
distributed independently in such a way that the probability $\bf
P\{n_i=1\}$ equals the given number $x_i\in[0,1],i=1,2,3$.

$L$-ensembles seems to be a more flexible means than Lefebvre's
real number function $f$ for some aspects. For example, let us
consider how "golden section" for categorization of stimuli
without measurable intensity can be explained in terms of
Lefebvre's theory (Lefebvre, 1995, p.51) and in terms of
$PL$-ensembles.

In this case Lefebvre adds equation $x_1=x_2,x_1=1-x_3$ to his
"Realist' condition" $x_3=f(x_1,x_2,x_3)$ (an justification is
given in (Lefebvre, 1995, p.51)). In turn, that yields the
equation $x_3^3-2x_3+1=0$ for the choice of $x_3$. One possible
solution is the well known "golden section" $x_3=\frac{\sqrt
5-1}{2}.$

Following the alternative approach, we construct the desired
$PL$-ensemble by first postulating the boolean "Realist'
condition" $n_3\to n_2\to n_1=n_3.$ Then considering the truth
area $R=\{000,001,010,101,111\}$ of the condition we form the
ensemble by means of boolean random variables $\zeta,\eta,\theta$
in the following fashion. The variables $\zeta,\eta$ are
independent with the mean value of $1-x_3$, and the value of the
random variable $\theta$ depending on the values of $\zeta,\eta$
as illustrated in Table \ref{table:1}.

\begin{table}
\begin{center}
\begin{tabular}{|c|c|c|}
\multicolumn{3}{r}{}\\
\hline
$\zeta$ & $\eta$ & $\theta$ \\
\hline
0 & 0 & 0,1\\
\hline
0 & 1 & 0 \\
\hline
1 & 0 & 1  \\
\hline
1 & 1 & 1 \\
\hline
\end{tabular}
\end{center}
\label{table:1} \caption{The solution list for boolean equation
$n_3\to n_2\to n_1=n_3.$}
\end{table}

It is important that in the first line of the table value $1$ is
chosen with the probability of $x_3$. Thus if $x_3$ satisfies
$x_3^3-2x_3+1=0$ then we obtain the desired $PL$-ensemble. Indeed,
every element of the ensemble is a "Realist" and the probability
to encounter an $L$-Subject with parameters $(n_1,n_2,1)$ is
determined by solutions to the equation $x_3^3-2x_3+1=0.$ Finally,
we arrive at the "golden section" choosing the corresponding
solution exactly as it was done by Lefebvre.

We believe that the $L$-ensemble tool provides additional
opportunities for Lefebvre's theory and its applications. Indeed,
the ensemble structure is a powerful parameter for modelling
because it can vary even though the average values are fixed.

\subsection{Application of approximating forms}

Now we propose an alternative model of the binary choice for the
same inputs and outputs. This model is constructed within the
framework of the theory proposed in this paper.

First we need to represent every choice that a Subject makes in
Lefebvre's model as a solution of the corresponding
extremalization problem of the aforementioned kind. Second we will
provide an algorithm for extremalization that computs results in
concordance with Lefebvre's theory.

The problem to pose such a extremalization problem is not trivial.
However, in our case it can be solved easily on the basis of the
interpretation of $x_i,i=\overline{1,3},$ given by Lefebvre.
Indeed, on one hand $x_1,x_2,x_3$ are connected to the
motivations: $x_1$ corresponds to the impulse (we continue to use
Lefebvre's terms) induced by the external world, $x_2$ corresponds
to the impulse induced by Subject's experience, and, finally,
$x_3$ corresponds to Subject's will. On the other hand, the values
of these variables describe objectives of the impulses. Thus, in
Lefebvre's model (boolean value of variable $x_i$ equals to 1(0))
if and only if (motivation $x_i$ pushes the Subject to the
positive(negative) pole).

\begin{figure}
\begin{center}
\setlength{\unitlength}{1mm}
\begin{picture}(50,30)
\put(6,5){$\bold{x_3}$} \put(6,15){$\bold{x_2}$}
\put(6,25){$\bold{x_1}$} \put(9,7){\vector(1,0){11}}
\put(9,17){\vector(1,0){11}}\put(9,27){\vector(1,0){11}}
\put(20,5){\framebox(10,25){$f$}} \put(30,17){\vector(1,0){11}}
\put(42,16){$\bold{z}$} %\put(20,0){Fig.1}
\put(9,7){\vector(1,0){11}}
\end{picture}
\end{center}
\caption{. Subject structure in Lefebvre's model.}\label{fig:1.}
\end{figure}
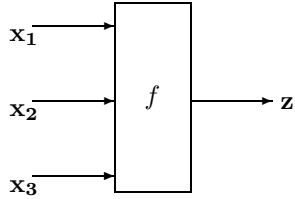

In order to avoid the ambiguity we will denote the boolean value
of variable $x_i$ in bold: $\bold{x_i}$. There are just eight
problems of choice in Lefebvre's model as:
$8=|\{(\bold{x_1,x_2,x_3})|\bold{x_i}\in\{0,1\},i=\overline{1,3}\}|.$
For each of these problems function $f$ computes a chosen pole
$\bold{z}$ (Figure \ref{fig:1.}).

It is easy to check that $\bold{z}\in\{\bold{x_1,x_2,x_3}\}$. So
in order to be accurate one needs to reconstruct impulse $z$ that
stands behind $\bold{z}$ and is implicit in Lefebvre's model. Thus
it is logical to think that Subject chooses one of the initial
impulse set which is defined by a current choice problem. Then the
Subject tries to implement the chosen impulse. Therefore within
the framework of the two-preference scheme we ought to set
$M=\{x_1,x_2,x_3\}$.

Then for a given boolean 3-tuple $(\bold{x_1,x_2,x_3})$ one needs
to propose a routine computing $\le_M,(L,\le_L)$, evaluating
mapping $\psi$, and an optimizing algorithm $\frak B$ in such a
way that for every input boolean 3-tuple algorithm $\frak B$
computes an extremum $x_i$ obeying condition
$\bold{x_i}=f(\bold{x_1,x_2,x_3})$. We can do this so that $\psi$
depends on $(\bold{x_1,x_2,x_3})$ only.

Our further consideration is based mainly on theorem
\ref{thm.thm2} which postulates the existence of a universal
representation of all mapping of kinds $\psi:M\to L$ when
preferences $(M,\le_M),(L,\le_L)$ are fixed. The representation
operates with the set of $\theta$-functions.

First, using the universality it is possible to define {\it\bf
any} evaluating mapping $\psi$ by means of an appropriate
substitution of $\theta$-functions into the corresponding
universal form. For that it is enough to link any
$(x_i,\bold{x_i})$ with an appropriate $\theta$- function.

Second, one needs to use the entire set $\Theta$. Taking into
account that different 3-tuples define different choice problems
we come to
$$
|\Theta|=|\{(x_i,\bold{x_i})|i=\overline{1,3},\bold{x_i}
\in\{0,1\}\}|=6.
$$
It is easy to see that this is possible only when
$(M,\le_M),(L,\le_L)$ are linear orderings and $|L|=2$. So we may
define the external preference by equalities:
$$
L=\{0,1\},\le_L=\{(0,0),(0,1),(1,1)\}.
$$

\begin{figure}
\begin{center}
\setlength{\unitlength}{1mm}
\begin{picture}(140,30)
\put(5,5){$\bold{x_3}$} \put(5,15){$\bold{x_2}$}
\put(5,25){$\bold{x_1}$} \put(11,7){\vector(1,0){5}}
\put(11,17){\vector(1,0){5}}\put(11,27){\vector(1,0){5}}
\put(18,5){$\theta_3^{\bold{x_3}}$}
\put(18,15){$\theta_2^{\bold{x_2}}$}
\put(18,25){$\theta_1^{\bold{x_1}}$} \put(25,7){\vector(1,0){5}}
\put(25,17){\vector(1,0){5}}\put(25,27){\vector(1,0){5}}
\put(31,5){\framebox(35,25){{\shortstack{"mixing" by\\
$\theta_1^{\bold{x_1}}\boxminus(\theta_2^{\bold{x_2}}\boxminus
\theta_3^{\bold{x_3}})$}}}} \put(66,17){\vector(1,0){5}}
\put(72,16){$\psi_{\bf{x_1,x_2,x_3}}$}
\put(87,17){\vector(1,0){5}}\put(92,5){\framebox(10,25){ $\frak
B$}} \put(102,17){\vector(1,0){5}} \put(108,16){$z$}
\put(110,17){\vector(1,0){5}} \put(116,16){$\bold{z}$}
\end{picture}
\end{center}
\caption{. Subject structure suggested in this paper.}
\label{fig:2}
\end{figure}
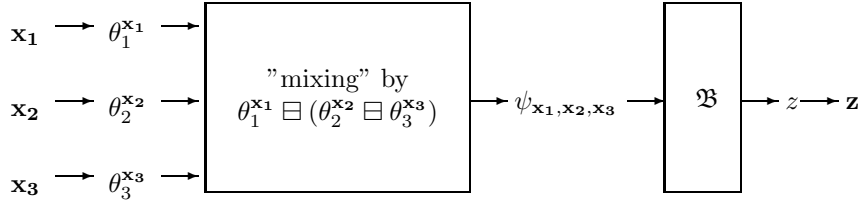

Further we choose the following linear order as the internal
preference:
$$
\le_M=\{(x_1,x_1),(x_2,x_2),(x_3,x_3),(x_1,x_2),(x_2,x_3),(x_1,x_3)\}.
$$
This is because in the considered case the internal (Subjective)
preference might be based on a degree of dependence of the states
on Subject's will. The world pressure $x_1$ depends on Subject to
the least extent. On the contrary, the dependence of $x_3$ on the
Subject is maximum. So the degree of dependence of $x_2$ on the
Subject lies in between the those two. The model can be now
finalized (Figure \ref{fig:2}).

Every choice problem the Subject is faced with can be
characterized by a certain boolean 3-tuple $(\bf{x_1,x_2,x_3})$ of
values of variables $x_1,x_2,x_3$. We associate a particular
evaluating mapping ("pure evaluation")
$\theta_i^{\bf{x_i}}\in\Theta$ with pair
$(x_i,\bold{x_i}),x_i,i=\overline{1,2,3},\bold{x_i}\in\{0,1\}.$
Thus:
\begin{center}
$\theta_i^{\bold{x_i}}(x_k) =
\begin{cases}
1,&  \text{ if }i<k,\\
\bold{x_i},& \text{ if } i=k,\\
0,&   \text{ if } k<i.
\end{cases}$
\end{center}

\noindent For a given external preference we can set
$\boxminus=\to^{\star}$ ($\to^{\star}$ is the boolean operation
dual to implication, see lemma 1 above). Hence every evaluation
mapping can be determined by the following formula:
$$
\psi_{\bold{x_1,x_2,x_3}}=
\theta_1^{\bold{x_1}}\to^{\star}(\theta_2^{\bold{x_2}}\to^{\star}
\theta_3^{\bold{x_3}}).
$$

\noindent Thus, one can consider evaluation mapping to be some
sort of "mixture" of pure evaluations. (It is worth noting a vague
analogy with quantum mechanics here. Lefebvre discussed a relation
of his model to the mechanics in (Lefebvre, 1991)).

Furthermore, $\psi_{(\bold{x_1,x_2,x_3})}:M_*\to\{0,1\}$
determines Subject's choice of state $x_i,i=1,2,3,$ for any given
choice problem and thereby the pole
$\psi_{(\bold{x_1,x_2,x_3})}(x_i)$.

In order to find a maximum value of a mapping of kind
$\psi':(M',\le_{M'})\to(L',\le_{L'})$ when
$(M',\le_{M'}),(L',\le_{L'})$ are linear orders, one can use an
easy algorithm based on the representation from theorem 2. We do
not formulate the algorithm or prove its correctness here.
Instead, we formulate a simplification of the algorithm $\frak A$
for $\psi_{\bold{x_1,x_2,x_3}}$ and $(M,\le_M),(L,\le_L)$ that we
have set above:

{\bf
\begin{enumerate}
\item Starting in the state $x_1$ in order $(M,\le_M)$
proceed to the nearest maximum $x_{j_1}$ of $\theta_1^{\bf{x_1}}$.
\item Continue from the state to the nearest
minimum $x_{j_2}$ of function $\theta_2^{\bf{x_2}}$ (due to its
place in the approximative form for $\psi_{\bold{x_1,x_2,x_3}}$).
\item Finally, starting from $x_{j_2}$ proceed to the nearest maximum $z$ of
$\theta_3^{\bold{x_3}}$.
\end{enumerate}}

\noindent The algorithm computes element $z$ of
$\arg\max\psi_{\bold{x_1,x_2,x_3}}$. Having the solution we know
the pole $\bold{z}$ chosen by Subject for parameters
$\bold{x_1,x_2,x_3}$. The results are presented in Table
\ref{table:2}.

Function $F(\bold{x_1,x_2,x_3})=\bold{z}$ which is presented by
the last column of the table deviates from Lefebvre's
$f(\bold{x_1,x_2,x_3})$ at one point only: $(0,1,1)$. In order to
make the two coherent it is sufficient to replace algorithm $\frak
A$ with an approximate algorithm $\frak a$. One can obtain the
latter algorithm by means of replacing words "nearest maximum"
("nearest minimum") with "nearest extremum in the maximizing
direction" ("nearest extremum in the minimizing direction") in the
description of $\frak A$ above:

{\bf \begin{enumerate}
\item Starting at the state $x_1$ in order $(M,\le_M)$ proceed to
the nearest extremum $x_{j_1}$ of $\theta_1^{\bold{x_1}}$ in the
maximizing direction.
\item Then continue from the state to the nearest
extremum $x_{j_2}$ of function $\theta_2^{\bold{x_2}}$ in the
minimizing direction (due to its place in the approximative form
for $\psi_{\bold{x_1,x_2,x_3}}$).
\item Finally, starting from $x_{j_2}$ proceed to the nearest extremum
$z$ of $\theta_3^{\bold{x_3}}$ in the maximizing direction.
\end{enumerate}}

\begin{table}
\begin{center}
\begin{tabular}{|c|c|c|c|c|c|c|}
\multicolumn{7}{r}{}\\
\hline \multicolumn{3}{|c}{\rule[-3mm]{0mm}{8mm}\bfseries
parameters} & \multicolumn{2}{|c|}{\rule[-3mm]{0mm}{8mm}\bfseries
choice by $\frak A $} &
\multicolumn{2}{|c|}{\rule[-3mm]{0mm}{8mm}\bfseries
choice by $\frak a$}\\
\hline $\bold{x_1}$ &$\bold{x_2}$ &$\bold{x_3}$ &$z$
&$F$ & $z$ & $f$ \\
\hline\hline
0 & 0 & 0 & $x_2$ & 0 & $x_2 $ & 0\\
\hline
0 & 0 & 1 & $x_3$ & 1 & $x_3 $ & 1\\
\hline
0 & 1 & 0 & $x_1$ & 0 & $x_1 $ & 0\\
\hline\hline
$\bold0$ & $\bold1$ & $\bold1$ & $x_3$ & $\bold1$ & $x_1 $ & $\bold0$\\
\hline\hline
1 & 0 & 0 & $x_1$ & 1 & $x_1 $ & 1\\
\hline
1 & 0 & 1 & $x_3$ & 1 & $x_1 $ & 1\\
\hline
1 & 1 & 0 & $x_1$ & 1 & $x_1 $ & 1\\
\hline
1 & 1 & 1 & $x_3$ & 1 & $x_1 $ & 1\\
\hline\hline
\end{tabular}
\end{center}
\caption{\ Results produced by algorithms $\frak A$ and $\frak a
$.}\label{table:2}
\end{table}

\noindent It turns out that boolean value $f(\bf{x_1,x_2,x_3})$
computed with algorithm $\frak a$ for all boolean 3-tuples
$(\bf{x_1,x_2,x_3})$ coincides with the value given by formula
$(x_3\to x_2)\to x_1$. Thus, algorithm $\frak a$ de facto
optimizes the external preference  in concordance with Lefebvre's
axioms. Indeed, if $\bold{x_1}=1$ then nothing happens: the start
state $x_1$ is the result of the choice. Hence, it is in
accordance with axiom $\mathcal L_3$. Otherwise, if
$\bold{x_1}=\bold{x_2}=0$ then $x_3$ or $x_2$ are chosen. In both
of these cases the boolean value of chosen variable coincides with
the $\bold{x_3}$ (axiom $\mathcal L_1$). Otherwise, $\bold{x_1}=0\
\&\ \bold{x_2}=1$ and the algorithm chooses $x_1$. This
corresponds to axiom $\mathcal L_2$. Thus, we are able to derive
these axioms from the algorithm.

\subsection{Discussion}

As shown above, the formula of human behavior proposed by Lefebvre
can be derived from our model given certain specific preferences
and optimization algorithm $\frak a$. Therefore, Lefebvre's
subjects appear distinguished merely by particular internal and
external orders $(M,\le_M),(L,\le_L)$.

Instead of evaluation mapping (pure evaluation) one may use
preference (pure preference) induced by it. In our model any
initial Subject's impulse $(x_i,\bf{x_i})$ is linked to the
partial order induced by mapping $\theta_{x_i}^{\bold{x_i}}$. This
order contributes to the external preference induced by
$\psi_{\bold{x_1,x_2,x_3}}$.

It should be noted that the statements are worded using 'extremes'
and not 'maxima' and 'minima'. This is so because the Subject can
use an approximation to the exact algorithm if the latter is
overly complex for it. Often such an approximation is sufficient
in practice.

One of the key strengths of our approach is the natural
generalization of the model for more than three states. In
particular, this is applicable when the Subject has two or more
levels of reflections.

Then, one can see that at the level of intentions (unlike the
level of their boolean values) there is a difference between the
case of $\bold{x_1}=\bold{x_2}=\bold{x_3}=0$ when the algorithm
$\frak a$ computes $x_2$ and the case of
$\bold{x_1}=\bold{x_2}=0,\bold{x_3}=1$ when the algorithm computes
$x_3$. Thus, it appears that we can't exactly follow Lefebvre's
reasoning on the "free will" when $\bold{x_1}=\bold{x_2}=0$.

If we adopt behavior function $F$ instead of $f$ then we would
lose the opportunity to explain the "golden section" effect
considered in the previous subsections. Therefore, in our model
the {\em inexact} algorithm appears to be the real cause of the
effect.

It may seem that $x_2<x_1$ ought to hold since we interpret $x_2$
as the "past experience" and $x_1$ as the "current pressure of the
environment". However, one should keep in mind that we are
currently dealing with an {\bf internal order} on states in the
process of decision-making. In that process "past experience"
$x_2$ serves the role of Subject's "current base" and it is $x_1$
that initiates decision-making. Variable $x_3$ is a means to
produce a solution and as such is most likely related to the
future.

\section{Conclusions}

As the paper demonstrates, the classical two-valued propositional
logic can be viewed as a realization of the method of successive
approximations for a decision-making problem within a framework of
Subject-in-an-environment survival.

Consequently, the classical propositional logic can take its
beginning from the survival problem. It is important that such
hypothetical origin of the logic appears quite natural.

Furthermore, this approach can serve as a background for
considering other families of mappings from one poset to another
with a chosen notion of simplicity of mapping. These families can
generate corresponding logics. So one may say that the
psychological effects described via Lefebvre's model considered
above can be interpreted as a logic rooted in evaluation functions
of the kind:
$$
\psi:(\{1,2,3\},\{1<2,1<3,2<3\})\to(\{0,1\},\{0<1\})
$$
implemented with a limited algorithm of extremum finding.


\begin{thebibliography}{99}


\bibitem{Birkhoff1967} Birkhoff, G. (1967),
"Lattice Theory", Providence/Rhode Island.

\bibitem{Bulitko2000}
Bulitko, V.K. (2000), Possible Origin of Logic, \textit{LANL},
arXiv:math.LO/0005050, 2000.

\bibitem{Bulitko97} Bulitko, V.K. (1997),
Lefebvre's Principle of Freedom and One Alternative Approach,
\textit{PSYCOLOQUY} \textbf{8(05)} human-choice.8.bulitko.

\bibitem{Kleene1967} Kleene, S.C. (1967)
"Mathematical Logic", John Wiley \& Sons Inc., New
York/London/Sydney.

\bibitem{Lefebvre1991} Lefebvre, V.A. (1991),
"The Formula of Man", Progress,(in Russian).

\bibitem{Lefebvre1995} Lefebvre, V.A. (1995),
The Anthropic Principle in Psychology and Human Choice,
\textit{PSYCOLOQUY} \textbf{6(29)}, human-choice.1.lefebvre.

\bibitem{Piaget1956}
Piaget, J. (1956), "Logic and psychology", Manchester University
Press.

\end{thebibliography}
\end{document}